 \documentclass[final,3p,times]{elsarticle}

\usepackage{amssymb}

\journal{Elsevier}

\usepackage{amsfonts,amsmath,amssymb,amsthm,stmaryrd,hyperref}
\usepackage[english]{babel}
\usepackage[T1]{fontenc}
\usepackage[ansinew]{inputenc}  
\usepackage{graphics}
\usepackage{graphicx}
\usepackage{tikz}
\usepackage{multicol}
\usetikzlibrary{patterns}
\usepackage{amssymb,subfigure,lineno,url}
\usepackage{float}

\newtheorem{theorem}{Theorem}

\newtheorem{lemma}{Lemma}
\newtheorem{proposition}{Proposition}

\newtheorem{definition}{Definition}
\newtheorem{remark}{Remark}

\begin{document}

\begin{frontmatter}



\title{Nonexistence of global weak solutions to semilinear wave equations involving time-dependent structural damping terms}

\author[ku,naam,cor1]{Mokhtar Kirane}
 \cortext[cor1]{Corresponding author}
 \ead{mokhtar.kirane@ku.ac.ae}
\affiliation[ku]{organization={Department of Mathematics},
            addressline={College of Arts and Sciences, Khalifa University of Science and Technology},
            city={Abu Dhabi},
            country={United Arab Emirates}}


\author[aum]{Ahmad Z. Fino}
\ead{ahmad.fino@aum.edu.kw}
\affiliation[aum]{organization={College of Engineering and Technology},
            addressline={American University of the Middle East}, 
            city={Egaila 54200},
            country={Kuwait}}

\author[squ]{Sebti Kerbal}
\ead{skerbal@squ.edu.om}

\affiliation[squ]{organization={Department of Mathematics},
            addressline={Sultan Qaboos University,  FracDiff Research Group (DR/RG/03)}, 
            city={Al-Khoud 123, Muscat},
            postcode={P.O. Box 36}, 
            country={Oman}}

\author[ku]{Aymen Laadhari}
 \ead{aymen.laadhari@ku.ac.ae}


\begin{abstract}
We consider a semilinear wave equation involving a time-dependent structural damping term of the form $\displaystyle\frac{1}{{(1+t)}^{\beta}}(-\Delta)^{\sigma/2} u_t$. Our results show the influence of the parameters $\beta,\sigma$ on the nonexistence of global weak solutions under assumptions on the given system data.
\end{abstract}




\end{frontmatter}


\section{Introduction}
The aim of the paper is to establish a nonexistence of global weak solutions  to the Cauchy
 problem for the following semilinear structurally damped wave equation
\begin{equation}\label{1}
\left\{\begin{array}{ll} \,\, \displaystyle {u_{tt}-\Delta
u +\frac{b_0}{(1+t)^{\beta}}(-\Delta)^{\sigma/2} u_t=|u|^p} &\displaystyle{x\in {\mathbb{R}^n},\,t>0,}\\
{}\\
\displaystyle{u(0,x)= u_0(x),\;\;u_t(0,x)=
u_1(x),\qquad\qquad}&\displaystyle{x\in {\mathbb{R}^n},}
\end{array}
\right.
\end{equation}
where $0<\sigma<2$, $n\geq 1$, $p>1$, $b_0$  is a positive constant, and $\beta\in\mathbb{R}$. Without loss of generality, we assume that $b_0=1$.\\

\medskip
Prior to our main blow-up results, we would like to emphasize that the special case of Eq. \eqref{1} for $b_0=0$ is governed by the classical semilinear wave equation, where the Strauss conjecture for this case states that: if $p \le p_S$ then there is no global solution for \eqref{1} under suitable sign assumptions on the initial data,  and for $p > p_S$ a global solution exists for small initial data; see e.g. \cite{John2,Strauss,YZ06,Zhou} among many other references, where $p_S$ is the positive  solution of the following quadratic equation
$(n-1)p^2-(n+1)p-2=0,$
and is given  by
$$
p_S=p_S(n):=\frac{n+1+\sqrt{n^2+10n-7}}{2(n-1)}.
$$
In order to have more general initial data, but still with compact support, Kato \cite{Kato} obtained with a slightly less sharp blow-up result an exponent of the form ${(n+1)}/{(n-1)}$ which is less than the Strauss exponent $p_S$, for $n\geq2$.\\
 We would like to emphasize here that the test function method, introduced by \cite{Zhang} and used by \cite{FinoKarch, Finokirane, PM}, plays a similar role as of Kato's method in the proof of blow-up results. In fact, the test function is effective in the case of parabolic equations as it provide exactly the critical exponent $p_c$, but in the case of hyperbolic equations (cf. \cite{PM}) we get the so-called Kato's exponent $p^*$ which is less than $p_c$. This is one of the weakness of the test function method but in general it can be applied to a more general equation and system with no positive condition on solutions.\\

When $\beta=0$, $\sigma\rightarrow2$, and $b_0=1$, problem \eqref{1} is reduced to
\begin{equation}\label{strongdamping}
 \displaystyle {u_{tt}-\Delta
u -\Delta u_t=|u|^p}, \ \ \ \ \displaystyle {x\in {\mathbb{R}^n},t>0,}
\end{equation}
which is called the viscoelastic damping case. D'Ambrosio and Lucente \cite[Theorem~4.2]{Ambrosio} proved that the solution of \eqref{strongdamping} blows-up in finite time when $1<p\leq (n+1)/(n-1)_+$, where $(\cdotp)_+:= \max\{0,\cdotp\}$, by applying the test function method. Similar result has been obtained recently  by Fino \cite{Fino} in the case of an exterior domain. On the other hand, D'Abbicco-Reissig \cite{Dabbicco} proved that there exists a global solution for \eqref{strongdamping} when $p>1+\frac{3}{n-1}$ $(n\geq 2)$ for sufficiently small initial data. Therefore, the exact value of the critical exponent is still an open question.
\\

Recently, problem \eqref{1} with $\sigma\rightarrow2$,
\begin{equation}\label{var}
 \displaystyle {u_{tt}-\Delta
u -\frac{b_0}{(1+t)^{\beta}}\Delta u_t=|u|^p}, \ \ \ \ \displaystyle {x\in {\mathbb{R}^n},t>0,}
\end{equation}
has been studied by Fino \& Hamza \cite{FH}. They proved that if $\int_{\mathbb{R}^n} u_1(x)dx > 0$,
$$
\left\{\begin{array}{ll}
p\in(1,\infty)&\hbox{for}\,\, n=1,2,\\
{}\\
 p\in(1, \frac{n}{n-2}]&\hbox{for}\,\, n\geq 3,\\
\end{array}
\right.
$$
 and
$$
\left\{\begin{array}{ll}
 1<p\leq \frac{n+1}{(n-1)_+}&\text{if}\,~ \beta\geq-1,\\
 &\\
1<p\leq \frac{n(1-\beta)+2}{(n(1-\beta)-2)_+}&\text{if}\,~ \beta\leq-1,\\
\end{array}
\right.
$$
then the mild solution of \eqref{var} blows-up in finite time.\\

The goal of this paper is to prove the nonexistence of global weak solutions of \eqref{1} which is a generalization of the results of \cite{FH} to the case of fractional Laplacian damping term by using a recent version of the test function method that relies on Lemma \ref{lemma3} below. The novelty of this paper is the choice of this test function which is consistent with the suppression of positivity usually assumed by many authors.\\

In the case of $\beta\neq0$ and $b_0=1$, we give an intuitive observation for understanding the influence of the damping  term ($\displaystyle\frac{1}{{(1+t)}^{\beta}}(-\Delta)^{\sigma/2} u_t$) by scaling argument. Let $u(t,x)$ be a solution of the linear strong damped wave equation
\begin{equation}\label{linear}
u_{tt}(t,x)-\Delta
u(t,x) +\frac{1}{(1+t)^{\beta}}(-\Delta)^{\sigma/2}  u_t(t,x)=0.
\end{equation}
When $\sigma\geq1-\beta$, we put
\begin{equation}\label{transf}
u(t,x)=v(\lambda(1+ t),\lambda x), \qquad \lambda (1+t)=s,\,\,\lambda x=y,
\end{equation}
with a parameter $\lambda>0$, we have
$$v_{ss}(s,y)-\Delta
v(s,y)+\frac{\lambda^{\sigma+\beta-1}}{s^{\beta}}(-\Delta)^{\sigma/2}  v_s(s,y)=0.$$
Thus, when $\sigma=1-\beta$ we notice that Eq. \eqref{linear} is invariant, while when $\sigma>1-\beta$, letting $\lambda\rightarrow 0$, we obtain the wave equation without damping
$$v_{ss}(s,y)-\Delta v(s,y)=0.$$
We note that $\lambda\rightarrow 0$ is corresponding to $t\rightarrow+\infty$.\\
On the other hand, when $\sigma<1-\beta$, we put
$$u(t,x)=v(\lambda^{\frac{\sigma}{1-\beta}}(1+ t),\lambda x), \qquad \lambda^{\frac{2}{1-\beta}}(1+t)=s,\,\,\lambda x=y,$$
with a parameter $\lambda>0$, we have
$$v_{ss}(s,y)+\frac{1}{s^{\beta}}(-\Delta)^{\sigma/2} v_s(s,y)-\lambda^{\frac{2(1-\beta-\sigma)}{1-\beta}}\Delta
v(s,y)=0.$$
In this case, letting $\lambda\rightarrow 0$, we obtain the pseudo-parabolic equation
$$v_{ss}(s,y)+\frac{1}{s^{\beta}}(-\Delta)^{\sigma/2} v_s(s,y)=0.$$
This means that the asymptotic behaviour of solutions depends on the behaviour of the coefficient of the damping term.\\

This paper is organized as follows. We start  in Sec. \ref{sec-main} by stating the main theorem of our work and we prove it in Sec. \ref{bmain}. Sec. \ref{preliminaries} is to collect some preliminaries.
\section{Main results}\label{sec-main}
\par
This section is aimed to state our main results. 

\begin{theorem}[\textbf{Blow-up}] \label{Blow-up}
 We assume that
$$(u_0,u_1)\in \big(L^1(\mathbb{R}^n)\cap H^1(\mathbb{R}^n)\big) \times \big(L^1(\mathbb{R}^n)\cap L^2(\mathbb{R}^n)\big) $$
satisfying the following condition:
\begin{equation} \label{initialdata}
\int_{\mathbb{R}^n} u_1(x)dx > 0.
\end{equation}
If 
\begin{equation}\label{30oc}
\left\{\begin{array}{ll}\displaystyle
 1<p<\min\left\{\frac{n+1}{(n-1)_+};\frac{n}{(n-\sigma)_+}\right\}&\text{if}\,~ \sigma\geq 1-\beta,\\
 &\\
 \displaystyle
1<p< \min\left\{\frac{n(1-\beta)+\sigma}{(n(1-\beta)-\sigma)_+};\frac{n}{(n-\sigma)_+}\right\}&\text{if}\,~ \sigma\leq 1-\beta,\\
\end{array}
\right.
\end{equation}
or
$$
\left\{\begin{array}{ll}
\displaystyle
 p=\frac{n+1}{n-1},&\displaystyle\quad\text{if}\,~ \sigma\geq 1-\beta,\,\,\sigma>\frac{2n}{n+1},\,\,n\geq2,\\
 &\\
 \displaystyle
p= \frac{n(1-\beta)+\sigma}{n(1-\beta)-\sigma},&\displaystyle\quad\text{if}\,~ n(1+\beta)<\sigma\leq 1-\beta,\,\,n\geq1,\\
\end{array}
\right.
$$
where $(\cdotp)_+:= \max\{0,\cdotp\}$, then problem \eqref{1} has no global weak solutions.
\end{theorem}

\begin{remark}
We note that, by taking the limit case $\sigma\rightarrow2$ we recover the same results of \cite{FH}.
\end{remark}

\begin{remark}\label{rk1bis}
We   stress that 
  the exponent  $\displaystyle\frac{n+1}{(n-1)_+}$ appearing  in \eqref{30oc} was  introduced  first in \cite{Kato} 
  to prove the nonexistence of global
solutions to the semilinear wave equation with the nonlinearity $|u|^p,$ subject to small initial data with compact support. 
\end{remark}

\section{Preliminaries}\label{preliminaries}
\par
\begin{definition}\cite{Silvestre}\label{def1}${}$\\
Let $\mathcal{S}$  be the Schwartz space of rapidly decaying $C^\infty$ functions in $\mathbb{R}^n$ and $s \in (0,1)$. The fractional Laplacian $(-\Delta)^s$ in $\mathbb{R}^n$ is a non-local operator defined on $\mathcal{S}$ by
\begin{eqnarray*}\displaystyle
 (-\Delta)^s v(x)&:=& \displaystyle C_{n,s}\,\, p.v.\int_{\mathbb{R}^n}\frac{v(x)- v(y)}{|x-y|^{n+2s}}\,dy\\\\
&=&\left\{\begin{array}{ll}
\displaystyle C_{n,s}\,\int_{\mathbb{R}^n}\frac{v(x)- v(y)}{|x-y|^{n+2s}}\,dy,&\quad\hbox{if}\,\,0<s<1/2,\\
{}\\
\displaystyle C_{n,s}\,\int_{\mathbb{R}^n}\frac{v(x)- v(y)-\nabla v(x)\cdotp(x-y)\mathcal{X}_{|x-y|<\delta}(y)}{|x-y|^{n+2s}}\,dy,\quad\forall\,\delta>0,&\quad\hbox{if}\,\,1/2\leq s<1,\\
\end{array}
\right.
\end{eqnarray*}
where $p.v.$ stands for Cauchy's principal value, and $\displaystyle C_{n,s}:= \frac{s\,4^s \Gamma(\frac{n}{2}+s)}{\pi^{\frac{n}{2}}\Gamma(1-s)}$.
\end{definition}

In fact, we are rarely going to use the fractional Laplacian operator in the Schwartz space; it can be extended to less regular functions as follows: for $s \in (0,1)$, $\varepsilon>0$, let
\begin{eqnarray*}
 \mathcal{L}_{s,\varepsilon}(\Omega)&:=&\left\{\begin{array}{ll}
\displaystyle L_s(\mathbb{R}^n)\cap C^{0,2s+\varepsilon}(\Omega)&\quad\hbox{if}\,\,0<s<1/2,\\
{}\\
\displaystyle L_s(\mathbb{R}^n)\cap C^{1,2s+\varepsilon-1}(\Omega),&\quad\hbox{if}\,\,1/2\leq s<1,\\
\end{array}
\right.
\end{eqnarray*}
where $\Omega$ be an open subset of $\mathbb{R}^n$, $C^{0,2s+\varepsilon}(\Omega)$ is the space of $2s+\varepsilon$- H\"{o}lder continuous functions on $\Omega$,  $C^{1,2s+\varepsilon-1}(\Omega)$ the space of functions of $C^1(\Omega)$ whose first partial
derivatives are H\"{o}lder continuous with exponent $2s+\varepsilon-1$, and
$$L_s(\mathbb{R}^n)=\left\{u:\mathbb{R}^n\rightarrow\mathbb{R}\quad\hbox{such that}\quad \int_{\mathbb{R}^n}\frac{u(x)}{1+|x|^{n+2s}}\,dx<\infty\right\}.$$

\begin{proposition}\label{Frac}\cite[Proposition~2.4]{Silvestre}${}$\\
Let $\Omega$ be an open subset of $\mathbb{R}^n$, $s \in (0,1)$, and $f\in \mathcal{L}_{s,\varepsilon}(\Omega)$ for some $\varepsilon>0$. Then $(-\Delta)^sf$ is a continuous function in $\Omega$ and $(-\Delta)^sf(x)$ is given by the pointwise formulas of Definition \ref{def1} for every $x\in\Omega$.
\end{proposition} 
\noindent{\bf Remark:} A simple sufficient condition for function $f$ to satisfy the conditions in Proposition \ref{Frac} is that $f\in  L^1_{loc}(\mathbb{R}^n)\cap C^{2}(\Omega)$.\\

Using \cite[Lemma~2.11]{DaoFino} and its proof, we have the following.
\begin{lemma}\label{lemma3}
Let $\langle x\rangle:=(1+(|x|-1)^4)^{1/4}$ for all $x\in\mathbb{R}^n$, $n\geq1$. Let $s \in (0,1]$ and $\phi:\mathbb{R}^n\rightarrow \mathbb{R}$ be a function defined by
\begin{equation}\label{testfunction}
\phi(x)=\left\{
\begin{array}{ll}
1&\hbox{if}\;\;|x|\leq1,\\
{}\\
\langle x\rangle^{-n-2s}&\hbox{if}\;\;|x|\geq1.
\end{array}
\right.
\end{equation}
Then
$$\phi\in \mathcal{C}^2(\mathbb{R}^n)\cap L^\infty(\mathbb{R}^n)\cap H^2(\mathbb{R}^n),\qquad \partial_x^2\phi\in L^\infty(\mathbb{R}^n),$$
 and the following estimate holds:
$$
\max\left\{\left|\Delta\phi(x)\right|,\,\left|(-\Delta)^s\phi(x)\right|\right\}\leq C\,\phi(x) \quad\hbox{for all}\;x\in\mathbb{R}^n.
$$
\end{lemma}

\begin{lemma}\label{lemma4}
Let $h$ be a smooth function satisfying $\partial_x^2h\in L^\infty(\mathbb{R}^n)$. For any $R>0$, let $h_R$ be a function defined by
$$ h_R(x):= h(R^{-1} x) \quad \text{ for all } x \in \mathbb{R}^n.$$
Then, $(-\Delta)^s (h_R)$,  $s \in (0,1]$,  satisfies the following scaling property
$$(-\Delta)^s h_R(x)= R^{-2s}(-\Delta)^sh(R^{-1} x), \quad \text{ for all } x \in \mathbb{R}^n. $$
\end{lemma}

\begin{lemma}\label{lemma5}\cite[Lemma~2.14]{DaoFino}\\
Let $s,\widetilde{s} \in (0,1]$, $R>0$ and $p>1$. Then, the following estimate holds
$$\int_{\mathbb{R}^n}(\phi_R(x))^{-\frac{1}{p-1}}\,\big|(-\Delta)^{\widetilde{s}} \phi_R(x)\big|^{\frac{p}{p-1}}\, dx\leq C\, R^{-\frac{2\widetilde{s}  p}{p-1}+n},$$
where $\phi_R(x):= \phi({x}/{R})$ and  $\phi$ as defined in (\ref{testfunction}).
\end{lemma}



\section{Proof of Theorem \ref{Blow-up}}\label{bmain}
The proof of Theorem \ref{Blow-up} relies mainly on the concept of weak solution of  the Cauchy problem (\ref{1}) and the use of the test function method. Let
$$X_{\delta,T}=\{\varphi\in C([0,\infty),H^2(\mathbb{R}^n))\cap C^1([0,\infty),H^\delta(\mathbb{R}^n))\cap C^2([0,\infty),L^2(\mathbb{R}^n)), \hbox{such that supp$\varphi\subset Q_T$}\},$$
where $Q_T:=[0,T]\times\mathbb{R}^n$, and the homogeneous fractional Sobolev space $H^\delta(\mathbb{R}^n)$, $\delta\in(0,2)$, is defined by
$$H^\delta(\mathbb{R}^n)=\left\{
\begin{array}{ll}
\{u\in L^2(\mathbb{R}^n); (-\Delta)^{\delta/2}u\in L^2(\mathbb{R}^n)\},&\,\,\hbox{if}\,\,\delta\in(0,1),\\\\
H^1(\mathbb{R}^n),&\,\,\hbox{if}\,\,\delta=1,\\\\
\{u\in H^1(\mathbb{R}^n); (-\Delta)^{\delta/2}u\in L^2(\mathbb{R}^n)\},&\,\,\hbox{if}\,\,\delta\in(1,2),\\\\
\end{array}
\right.$$
endowed with the norm
$$\|u\|_{H^\delta(\mathbb{R}^n)}=\left\{
\begin{array}{ll}
\|u\|_{L^2(\mathbb{R}^n)}+\left\|(-\Delta)^{\delta/2}u\right\|_{L^2(\mathbb{R}^n)},&\,\,\hbox{if}\,\,\delta\in(0,1),\\\\
\|u\|_{L^2(\mathbb{R}^n)}+\|\nabla u\|_{L^2(\mathbb{R}^n)},&\,\,\hbox{if}\,\,\delta=1,\\\\
\|u\|_{L^2(\mathbb{R}^n)}+\|\nabla u\|_{L^2(\mathbb{R}^n)}+\left\|(-\Delta)^{\delta/2}u\right\|_{L^2(\mathbb{R}^n)},&\,\,\hbox{if}\,\,\delta\in(1,2).\\\\
\end{array}
\right.$$
The weak formulation associated with (\ref{1}) reads  as follows:

\begin{definition}(Weak solution)\\
Let $T>0$, and $u_0,u_1\in L^2(\mathbb{R}^n)$. A function
$$u\in L^1((0,T),L^{2}(\mathbb{R}^n))\cap L^p((0,T);L^{2p}(\mathbb{R}^n)),$$
 is said to be a weak solution of $(\ref{1})$ on $[0,T)\times\mathbb{R}^n$ if 
 \begin{eqnarray*}
&{}&\int_{Q_T}|u|^p\varphi\,dt\,dx+\int_{\mathbb{R}^n}u_1(x)\varphi(0,x)\,dx\\
&{}&+\int_{\mathbb{R}^n}u_0(x)(-\Delta)^{\sigma/2}\varphi(0,x)\,dx-\int_{\mathbb{R}^n}u_0(x)\varphi_t(0,x)\,dx\\
&{}&=\int_{Q_T}u\varphi_{tt}\,dt\,dx-\int_0^T\frac{1}{(1+t)^{\beta}}\int_{\mathbb{R}^n}u\,(-\Delta)^{\sigma/2}\varphi_t\,dx\,dt\\
&{}&-\int_{Q_T}u \Delta \varphi\,dt\,dx+\int_0^T\frac{\beta}{(1+t)^{\beta+1}}\int_{\mathbb{R}^n}u\,(-\Delta)^{\sigma/2} \varphi\,dx\,dt,
\end{eqnarray*}
holds for all $\varphi \in X_{\sigma,T}$. We denote the lifespan for the weak solution by
$$T_w(u_0,u_1):=\sup\{T\in(0,\infty]\,\,\hbox{ for which there exists a unique weak solution u to (\ref{1})}\}.$$
Moreover, if $T>0$ can be arbitrary chosen, i.e. $T_w(u_0,u_1)=\infty$, then $u$ is called a global weak solution of (\ref{1}).
\end{definition}
\noindent{\bf Proof of Theorem \ref{Blow-up}.}
Let $u$ be a global weak solution of \eqref{1}, that is
  \begin{eqnarray}\label{weak1}
  &{}&\int_{Q_T}|u|^p\varphi\,dt\,dx+\int_{\mathbb{R}^n}u_1(x)\varphi(0,x)\,dx\nonumber\\
&{}&+\int_{\mathbb{R}^n}u_0(x)(-\Delta)^{\sigma/2}\varphi(0,x)\,dx-\int_{\mathbb{R}^n}u_0(x)\varphi_t(0,x)\,dx\nonumber\\
&{}&=\int_{Q_T}u\varphi_{tt}\,dt\,dx-\int_0^T\frac{1}{(1+t)^{\beta}}\int_{\mathbb{R}^n}u\,(-\Delta)^{\sigma/2}\varphi_t\,dx\,dt\nonumber\\
&{}&-\int_{Q_T}u \Delta \varphi\,dt\,dx+\int_0^T\frac{\beta}{(1+t)^{\beta+1}}\int_{\mathbb{R}^n}u\,(-\Delta)^{\sigma/2} \varphi\,dx\,dt,
\end{eqnarray}
for all $T>0$, and all $\varphi\in X_{\sigma,T}$.

Let $T>0$. Now, we introduce the following test function
\begin{equation}\label{psi}
\varphi(t,x)=\psi^\eta(t)\phi_{T^d}(x)=\psi^\eta(t)\phi\left(\frac{x}{T^d}\right)
\end{equation}
where $\phi$ is defined in (\ref{testfunction}), $\psi(t)=\Psi(\frac{t}{T})$, $\eta\gg1$, $d>0$ are constants that will be determined later, and $\Psi\in C^\infty(\mathbb{R}_+)$ is a cut-off non-increasing
function such that
$$\Psi(r)=\left\{\begin {array}{ll}\displaystyle{1}&\displaystyle{\quad\mbox{if }\,0\leq r\leq 1/2,}\\\\
\displaystyle{\searrow}&\displaystyle{\quad\mbox{if }\,1/2\leq r\leq 1,}\\\\
\displaystyle{0}&\displaystyle{\quad\mbox {if }\,r\geq 1.}
\end {array}\right.$$\\
From the formulation \eqref{weak1}, we get the following inequality
\begin{eqnarray}\label{3}
&{}&\int_{Q_T}|u|^p\varphi\,dt\,dx+\int_{\mathbb{R}^n}u_1(x)\varphi(0,x)\,dx\nonumber\\
&{}&\leq \int_{\frac{T}{2}}^T\int_{\mathbb{R}^n}|u|\,|\varphi_{tt}|\,dx\,dt+\int_{Q_T}|u|\,|\Delta \varphi|\,dt\,dx\nonumber\\
&{}&\,\,+\int_{\frac{T}{2}}^T\frac{1}{(1+t)^{\beta}}\int_{\mathbb{R}^n}|u||(-\Delta)^{\sigma/2}\varphi_t|\,dx\,dt+\int_0^T
\frac{\beta}{(1+t)^{\beta+1}}\int_{\mathbb{R}^n}|u|\,\,|(-\Delta)^{\sigma/2} \varphi|\,dx\,dt\nonumber\\
&{}&\,\,+\int_{\mathbb{R}^n}|u_0|\left(|(-\Delta)^{\sigma/2}\varphi(0,x)|+|\varphi_t(0,x)|\right)\,dx\nonumber\\
&{}&=:I_1+I_2+I_3+I_4+I_5.
\end{eqnarray}
Let $\varepsilon>0$. By applying $\varepsilon$-Young's inequality 
$$
AB\leq\varepsilon A^p+C(\varepsilon,p)B^{p^\prime},\,\, A\geq0,\;B\geq0,\;p+p^\prime=pp^\prime,\,\, C(\varepsilon,p)=(p-1)(\varepsilon\,p^p)^{-1/(p-1)},
$$
we obtain the estimation for the first integral
\begin{eqnarray}\label{4}
I_1&\leq&\int_{Q_T}|u|\varphi^{1/p}\varphi^{-1/p}\phi_{T^d}
|(\psi^\eta)_{tt}|\,dt\,dx\nonumber\\
&\leq&\varepsilon\int_{Q_T}|u|^{p}\varphi\,dt\,dx+C\int_{Q_T}\phi_{T^d}\psi^{-\eta/(p-1)}
|(\psi^\eta)_{tt}|^{p^\prime}\,dt\,dx.
\end{eqnarray} 
As $(\psi^\eta)_{tt}=\eta\psi^{\eta-1}\psi_{tt}+\eta(\eta-1)\psi^{\eta-2}|\psi_t|^2$, the inequality \eqref{4} becomes
\begin{equation}\label{4bis}
I_1\leq\varepsilon\int_{Q_T}|u|^{p}\psi\,dt\,dx+C\int_{Q_T}\phi_{T^d}\psi^{\eta-p^\prime}|\psi_{tt}|^{p^\prime}\,dt\,dx+\,C\int_{Q_T}\phi_{T^d}\psi^{\eta-2p^\prime}|\psi_{t}|^{2p\prime}\,dt\,dx.
\end{equation} 
Proceeding similarly as for \eqref{4bis},  we get
\begin{equation}\label{5}
I_2\leq\varepsilon \int_{Q_T}|u|^{p}\varphi\,dt\,dx+C\int_{Q_T}\psi^\eta\phi^{-1/(p-1)}_{T^d}
|\Delta\phi_{T^d}|^{p^\prime}\,dt\,dx.
\end{equation} 
In the same way, we write
\begin{eqnarray}\label{6}
I_3&=&\int_{\frac{T}{2}}^T
\frac{1}{(1+t)^{\beta}}\int_{\mathbb{R}^n}|u|\varphi^{1/p}\varphi^{-1/p}\,|(-\Delta)^{\sigma/2}\phi_{T^d}||(\psi^\eta)_t|\,dx\,dt\nonumber\\
&\leq&\varepsilon \int_{Q_T}|u|^{p}\varphi\,dt\,dx\nonumber\\
&{}&+\,C\int_{\frac{T}{2}}^T
\frac{1}{(1+t)^{\beta p^\prime}}\int_{\mathbb{R}^n}
\phi_{T^d}^{-1/(p-1)}
|(-\Delta)^{\sigma/2}\phi_{T^d}|^{p^\prime}\,\psi^{\eta-p^\prime}\,|\psi_t|^{p^\prime}\,dx\,dt.\quad
\end{eqnarray} 
 Clearly,
$$
\frac{1}{(1+t)^{\beta p^\prime}}\leq C\,T^{-\beta p^\prime},\qquad\forall\,t\in\left(\frac{T}{2},T\right),$$
 therefore, 
\begin{equation}\label{6bis}
I_3\leq\varepsilon \int_{Q_T}|u|^{p}\varphi\,dt\,dx+C\,T^{-\beta p^\prime}\, \int_{Q_T}\phi_{T^d}^{-1/(p-1)}
|(-\Delta)^{\sigma/2}\phi_{T^d}|^{p^\prime}\,\psi^{\eta-p^\prime}\,|\psi_t|^{p^\prime}\,dt\,dx.
\end{equation} 
In the same manner,
\begin{eqnarray}\label{7}
I_4&\leq&C\int_{0}^T\frac{1}{(1+t)^{\beta+1}}\int_{\mathbb{R}^n}|u|\varphi^{1/p}\varphi^{-1/p}\,\,\psi^\eta\,|(-\Delta)^{\sigma/2}\phi_{T^d}|\,dx\,dt\nonumber\\
&\leq&\varepsilon\int_{Q_T}|u|^{p}\varphi\,dt\,dx+C\,\int_0^T\frac{1}{(1+t)^{(\beta+1)p^\prime}}\int_{\mathbb{R}^n}
\phi_{T^d}^{-1/(p-1)}|(-\Delta)^{\sigma/2}\phi_{T^d}|^{p^\prime}\,\psi^{\eta}\,dx\,dt.
\end{eqnarray} 
Finally, it remains only to control the term $I_5$. By exploiting  the identity $\varphi_t (0,x)=\eta\psi_t(0) \phi_{T^d}(x)$, we infer
\begin{equation}\label{8}
I_5 \leq C\,\int_{\mathbb{R}^n}|u_0|\left(|(-\Delta)^{\sigma/2}\phi_{T^d}|+|\psi_t(0)|\phi_{T^d}\right)\,dx.
\end{equation} 
Plugging \eqref{3} together  with \eqref{4bis}-\eqref{8} and  choosing $\varepsilon$ small enough, we deduce  that
\begin{eqnarray*}
&{}&\int_{Q_T}|u|^p\varphi\,dx\,dt+\int_{\mathbb{R}^n}u_1(x)\phi_{T^d}(x)\,dx\\
&{}&\leq\,C\int_{Q_T}\phi_{T^d}\psi^{\eta-p^\prime}|\psi_{tt}|^{p^\prime}\,dt\,dx+\,C\int_{Q_T}\phi_{T^d}\psi^{\eta-2p^\prime}|\psi_{t}|^{2p\prime}\,dt\,dx\\
&{}&\quad+\,C\int_{Q_T}\psi^\eta\phi^{-1/(p-1)}_{T^d}
|\Delta\phi_{T^d}|^{p^\prime}\,dt\,dx\\
&{}&\quad+\,C\,T^{-\beta p^\prime}\, \int_{Q_T}\phi_{T^d}^{-1/(p-1)}
|(-\Delta)^{\sigma/2}\phi_{T^d}|^{p^\prime}\,\psi^{\eta-p^\prime}\,|\psi_t|^{p^\prime}\,dt\,dx\\
&{}&\quad+\,C\,\int_0^T\frac{1}{(1+t)^{(\beta+1)p^\prime}}\int_{\mathbb{R}^n}
\phi_{T^d}^{-1/(p-1)}|(-\Delta)^{\sigma/2}\phi_{T^d}|^{p^\prime}\,\psi^{\eta}\,dx\,dt\\
&{}&\quad+\,C\,\int_{\mathbb{R}^n}|u_0|\left(|(-\Delta)^{\sigma/2}\phi_{T^d}|+|\psi_t(0)|\phi_{T^d}\right)\,dx.
\end{eqnarray*}
Taking account of the expression of $\varphi$ given by \eqref{psi} and Lemmas \ref{lemma3}-\ref{lemma5},
we infer that
\begin{eqnarray}\label{10}
&{}&\int_{Q_T}|u|^p\varphi\,dx\,dt+\int_{\mathbb{R}^n}u_1(x)\phi_{T^d}(x)\,dx\nonumber\\
&{}&\leq C\;T^{-2p^\prime+1+nd}+C\;T^{-2dp^\prime+1+nd}+C\;T^{-\beta p^\prime-p^\prime-\sigma dp^\prime+1+nd}\\
&&\quad +\,C\;T^{-\sigma dp^\prime+nd}
\int_0^T(1+t)^{-\frac{(\beta+1)\,p}{p-1}}\,dt+\,C\left(T^{-\sigma d}+T^{-1}\right)\,\int_{\mathbb{R}^n}|u_0(x)|\,dx.\nonumber
\end{eqnarray}

Since $0<\sigma<2$, we notice that the cases of $\sigma\geq1-\beta$ and $\sigma\leq1-\beta$ are equivalent to 
$$\sigma\geq1-\beta,\quad\hbox{for all}\,\,\beta>-1, \quad\qquad\hbox{and}\quad\qquad \sigma<1-\beta,\quad\hbox{for all}\,\,\beta<1,$$
therefore, we distinguish two cases:\\

 \noindent {\bf I. Case: $\sigma\geq1-\beta$, for all $\beta>-1$.}\\
In this case, we choose $d=1$. \\
  
\noindent \underline{Subcritical case $\displaystyle p<\min\{\frac{n+1}{(n-1)_+};\frac{n}{(n-\sigma)_+}\}$.}\\
 Note that,
$$
\int_0^T(1+t)^{\displaystyle -\frac{(\beta+1)\,p}{p-1}}\,dt\leq C\,\left\{\begin{array}{ll}\displaystyle
T^{\displaystyle 1-\frac{(\beta+1)\,p}{p-1}}&\text{if}~\beta\,p<-1,\\
\ln T&\text{if}~\beta\,p=-1,\\
1&\text{if}~\beta\,p>-1.\\
\end{array}
\right.
$$
We have two cases to distinguish.
 \begin{enumerate}
 \item[a)] \underline{If $\beta\geq0$.} In this case, we have $\beta\,p>-1$ and so
 $$\int_0^T(1+t)^{\displaystyle -\frac{(\beta+1)\,p}{p-1}}\,dt\leq C.$$ 
 Therefore, (\ref{10}) implies
\begin{eqnarray*}
&{}&\int_{Q_T}|u|^p\varphi\,dx\,dt+\int_{\mathbb{R}^n}u_1(x)\phi_{T}(x)\,dx\\
&{}&\leq C\;T^{-2 p^\prime+1+n}+C\;T^{-(\beta+1) p^\prime-\sigma p^\prime+1+n}+ C\;T^{-\sigma p^\prime+n}+
\,C\left(T^{-\sigma}+T^{-1}\right)\,\int_{\mathbb{R}^n}|u_0(x)|\,dx.
\end{eqnarray*}
Using the fact $\beta\,p>-1\Longrightarrow{-(\beta+1) p^\prime-\sigma p^\prime+1+n}<{-\sigma p^\prime+n}$, we conclude that
\begin{eqnarray}\label{10d11}
&{}&\int_{Q_T}|u|^p\varphi\,dx\,dt+\int_{\mathbb{R}^n}u_1(x)\phi_{T}(x)\,dx\nonumber\\
&{}&\leq C\;T^{-2p^\prime+1+n}+ C\;T^{-\sigma p^\prime+n}+
\,C\left(T^{-\sigma}+T^{-1}\right)\,\int_{\mathbb{R}^n}|u_0(x)|\,dx.
\end{eqnarray}
Note that,  we can easily see that
$$-2p^\prime+1+n<0\Longleftrightarrow p<\frac{n+1}{(n-1)_+}\quad\hbox{and}\quad-\sigma p^\prime+n<0\Longleftrightarrow p<\frac{n}{(n-\sigma)_+}.$$
Letting $T\rightarrow\infty$, and using the Lebesgue
dominated convergence theorem together with $u_1\in L^1(\mathbb{R}^n)$, we conclude that
$$\int_{\mathbb{R}^n}u_1(x)\,dx\leq 0.$$
This contradicts our assumption \eqref{initialdata}.
 \item[b)] \underline{If $-1<\beta<0$.} We have three cases for $n$.
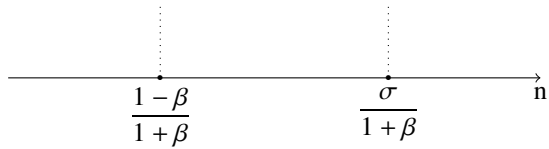
  \begin{center}
\begin{tikzpicture}[scale=1]
 \draw[->] (-4,0) -- (3,0) node[below]{n};
 \draw (1,0) node[below]{$\displaystyle\frac{\sigma}{1+\beta}$};
  \draw (1,0) node{{\tiny$\bullet$}};
   \draw[dotted](1,0)--(1,1);
    \draw (-2,0) node[below]{$\displaystyle\frac{1-\beta}{1+\beta}$};
  \draw (-2,0) node{{\tiny$\bullet$}};
   \draw[dotted](-2,0)--(-2,1);

        \end{tikzpicture}
 \end{center}
 Note that when $n=1$, we have $n=1<\frac{1-\beta}{1+\beta}$.
   \begin{enumerate}
 \item[i)] If $\displaystyle n\geq\frac{\sigma}{1+\beta}$, then $\displaystyle n\geq\frac{1-\beta}{1+\beta}$ and therefore
 $$p<\min\left\{\frac{n+1}{n-1};\frac{n}{n-\sigma}\right\}\leq-\frac{1}{\beta},\quad\hbox{for all}\,\,n\geq2,$$
 i.e. $\beta\,p>-1$, which implies a contradiction by following the same calculations as in part a).
 \item[ii)] If $\displaystyle\frac{1-\beta}{1+\beta}<n<\frac{\sigma}{1+\beta}$ when $\sigma>1-\beta$, then we have
  $$n<\frac{\sigma}{1+\beta}<\frac{\sigma}{2-\sigma},$$
 and therefore
  $$p<\min\left\{\frac{n+1}{n-1};\frac{n}{n-\sigma}\right\}=\frac{n+1}{n-1}<-\frac{1}{\beta},\quad\hbox{for all}\,\,n\geq2,$$
 i.e. $\beta\,p>-1$, which implies a contradiction by following the same calculations as in part a).
   \item[iii)] If $\displaystyle n\leq\frac{1-\beta}{1+\beta}$, then $\displaystyle n\leq\frac{\sigma}{1+\beta}$, i.e.
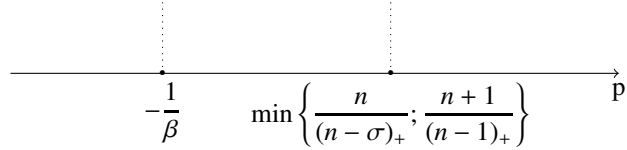
    \begin{center}
\begin{tikzpicture}[scale=1]
 \draw[->] (-4,0) -- (4,0) node[below]{p};
 \draw (1,0) node[below]{$\min\left\{\displaystyle\frac{n}{(n-\sigma)_+};\frac{n+1}{(n-1)_+}\right\}$};
  \draw (1,0) node{{\tiny$\bullet$}};
   \draw[dotted](1,0)--(1,1);
 \draw (-2,0) node[below]{$\displaystyle -\frac{1}{\beta}$};
  \draw (-2,0) node{{\tiny$\bullet$}};
   \draw[dotted](-2,0)--(-2,1);
        \end{tikzpicture}
 \end{center}

$\bullet$ If $p\leq -1/\beta$, i.e. $p\beta\geq -1$, then 
$$\int_0^T(1+t)^{\displaystyle -\frac{(\beta+1)\,p}{p-1}}\,dt\leq C\,\ln T,\qquad\hbox{for all}\,\, T\gg1.$$ 
Then, (\ref{10}) implies
\begin{eqnarray*}
&{}&\int_{Q_T}|u|^p\varphi\,dx\,dt+\int_{\mathbb{R}^n}u_1(x)\phi_{T}(x)\,dx\\
&{}&\leq C\;T^{-2 p^\prime+1+n}+C\;T^{-(\beta+1) p^\prime-\sigma p^\prime+1+n}+ C\;T^{-\sigma p^\prime+n}\,\ln T+
\,C\left(T^{-\sigma}+T^{-1}\right)\,\int_{\mathbb{R}^n}|u_0(x)|\,dx,
\end{eqnarray*}
 for all $T\gg1$. As $p\beta\geq -1\Longrightarrow (\beta+1)p^\prime\geq 1$, we conclude that
\begin{eqnarray}\label{10d12}
&{}&\int_{Q_T}|u|^p\varphi\,dx\,dt+\int_{\mathbb{R}^n}u_1(x)\phi_{T}(x)\,dx\nonumber\\
&{}&\leq C\;T^{-2p^\prime+1+n}+ C\;T^{-\sigma p^\prime+n}\,\ln T+
\,C\left(T^{-\sigma}+T^{-1}\right)\,\int_{\mathbb{R}^n}|u_0(x)|\,dx, \quad \forall\,T\gg1.
\end{eqnarray}
Note that,  we can easily see that
$$-2p^\prime+1+n<0\Longleftrightarrow p<\frac{n+1}{(n-1)_+}\quad\hbox{and}\quad-\sigma p^\prime+n<0\Longleftrightarrow p<\frac{n}{(n-\sigma)_+}.$$
Letting $T\rightarrow\infty$, using the fact that $\ln T\leq T^{\frac{\sigma p^\prime-n}{2}}$ (because $p<\frac{n}{(n-\sigma)_+}$) and the Lebesgue
dominated convergence theorem, we conclude that
$$\int_{\mathbb{R}^n}u_1(x)\,dx\leq 0.$$
This contradicts our assumption \eqref{initialdata}.\\

$\bullet$ If $p> -1/\beta$, i.e. $p\beta< -1$, then 
$$\int_0^T(1+t)^{-\frac{(\beta+1)\,p}{p-1}}\,dt\leq C\,T^{1-\frac{(\beta+1)\,p}{p-1}},\qquad\hbox{for all}\,\, T>0.$$ 
Then, (\ref{10}) implies
\begin{eqnarray*}
&{}&\int_{Q_T}|u|^p\varphi\,dx\,dt+\int_{\mathbb{R}^n}u_1(x)\phi_{T}(x)\,dx\\
&{}&\leq C\;T^{-2 p^\prime+1+n}+C\;T^{-(\beta+1) p^\prime-\sigma p^\prime+1+n}+
\,C\left(T^{-\sigma}+T^{-1}\right)\,\int_{\mathbb{R}^n}|u_0(x)|\,dx,
\end{eqnarray*}
 for all $T\gg1$. As $\sigma\geq1-\beta\Longrightarrow n\leq\frac{1-\beta}{1+\beta}\leq\frac{\sigma}{2-\sigma}$, we conclude that
 $$p<\min\left\{\frac{n+1}{(n-1)_+};\frac{n}{(n-\sigma)_+}\right\}=\frac{n+1}{(n-1)_+}\leq\frac{n+1}{(n-\beta-\sigma)_+}\Longrightarrow -(\beta+1) p^\prime-\sigma p^\prime+1+n<0,$$
and
$$p<\min\left\{\frac{n+1}{(n-1)_+};\frac{n}{(n-\sigma)_+}\right\}=\frac{n+1}{(n-1)_+}\Longleftrightarrow -2p^\prime+1+n<0.$$
Note that $\sigma\geq 1-\beta>1$, and $2>\sigma>\sigma+\beta$, therefore 
$$
\left\{\begin{array}{ll}
(n-\beta-\sigma)_+=(n-1)_+=(n-\sigma)_+=0,&\quad\hbox{when}\,\,n=1,\\\\
(n-\beta-\sigma)_+=n-\beta-\sigma,\,\,(n-1)_+=n-1,\,\,(n-\sigma)_+=n-\sigma,&\quad\hbox{when}\,\,n\geq2,\\
\end{array}
\right.
$$
 Letting $T\rightarrow\infty$, using the Lebesgue
dominated convergence theorem, we conclude that
$$\int_{\mathbb{R}^n}u_1(x)\,dx\leq 0.$$
This contradicts our assumption \eqref{initialdata}.\\
 
 \end{enumerate}
\end{enumerate}

 \noindent \underline{Critical case: $p<\infty$ and $\sigma\geq1$ when $n=1$, or $\displaystyle p=\frac{n+1}{n-1}$ and $\displaystyle\sigma>\frac{2n}{n+1}$ when $n\geq2$.}\\
Note that when $\displaystyle 2\leq n<\frac{\sigma}{2-\sigma}$, we have
 $$p=\frac{n+1}{n-1}<\frac{n}{n-\sigma}\quad\hbox{i.e.}\quad  -\sigma p^\prime+n<0,$$
  and when $\sigma\geq1$ and $n=1$, we also have $-\sigma p^\prime+n<0$.\\
We have two cases to distinguish.
 \begin{enumerate}
 \item[a)] \underline{If $\beta\geq0$.} In this case, we have $\beta\,p>-1$ and so
 $$\int_0^T(1+t)^{\displaystyle -\frac{(\beta+1)\,p}{p-1}}\,dt\leq C.$$ 
 From the subcritical case, we can see that we have
 \begin{equation}\label{regularity1}
 u\in L^p((0,\infty);L^p(\mathbb{R}^n)).
 \end{equation}
 On the other hand, by applying H\"{o}lder's inequality instead of Young's inequality, we get
\begin{eqnarray*}
\int_{\mathbb{R}^n}u_1(x)\phi_{T}(x)\,dx &\leq& C\left(\int_{\frac{T}{2}}^T\int_{\mathbb{R}^n}|u|^{p}\varphi\,dx\,dt\right)^{1/p}+\,C\left(\int_0^T\int_{|x|\geq T}|u|^{p}\varphi\,dx\,dt\right)^{1/p}\\
&{}&+ C\;T^{-\sigma p^\prime+n}+\,C\left(T^{-\sigma}+T^{-1}\right)\,\int_{\mathbb{R}^n}|u_0(x)|\,dx.
\end{eqnarray*}
 Letting $T\longrightarrow\infty$ and taking into consideration \eqref{regularity1}, we get a contradiction.
 \item[b)] \underline{If $-1<\beta<0$.} We have three cases for $n$.
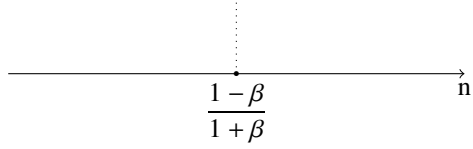
  \begin{center}
\begin{tikzpicture}[scale=1]
 \draw[->] (-3,0) -- (3,0) node[below]{n};
    \draw (0,0) node[below]{$\displaystyle\frac{1-\beta}{1+\beta}$};
  \draw (0,0) node{{\tiny$\bullet$}};
   \draw[dotted](0,0)--(0,1);

        \end{tikzpicture}
 \end{center}
 Note that when $n=1$, we have $\displaystyle n=1<\frac{1-\beta}{1+\beta}$.
  \begin{enumerate}
 \item[i)] If $\displaystyle n>\frac{1-\beta}{1+\beta}$, then
 $$p=\frac{n+1}{n-1}<-\frac{1}{\beta},$$
 i.e. $\beta\,p>-1$, which implies a contradiction by following the same calculations as in part a).
  \item[ii)] If $n=\frac{1-\beta}{1+\beta}$, then 
   $$p=\frac{n+1}{n-1}=-\frac{1}{\beta},$$
 i.e. $\beta\,p=-1$, and therefore
$$\int_0^T(1+t)^{\displaystyle -\frac{(\beta+1)\,p}{p-1}}\,dt\leq C\,\ln T,\qquad\hbox{for all}\,\, T>0.$$ 
 From the subcritical case \eqref{10d11}, we can see easily that we have
 \begin{equation}\label{regularity2}
 u\in L^p((0,\infty);L^p(\mathbb{R}^n)).
 \end{equation}
 On the other hand, by applying H\"{o}lder's inequality instead of Young's inequality and using $\displaystyle p=\frac{n+1}{n-1}$, we get
\begin{eqnarray*}
\int_{\mathbb{R}^n}u_1(x)\phi_{T}(x)\,dx &\leq& C\left(\int_{\frac{T}{2}}^T\int_{\mathbb{R}^n}|u|^{p}\varphi\,dx\,dt\right)^{1/p}+\,C\left(\int_0^T\int_{|x|\geq T}|u|^{p}\varphi\,dx\,dt\right)^{1/p}\\
&{}&+ C\;T^{-\sigma p^\prime+n}\ln T+\,C\left(T^{-\sigma}+T^{-1}\right)\,\int_{\mathbb{R}^n}|u_0(x)|\,dx.
\end{eqnarray*}
 Letting $T\longrightarrow\infty$ and taking into consideration \eqref{regularity2} and the fact that $\displaystyle \ln T\leq T^{\displaystyle\frac{\sigma p^\prime-n}{2}}$ (because $\displaystyle p<\frac{n}{n-\sigma}$), we get
$$\int_{\mathbb{R}^n}u_1(x)\,dx\leq 0.$$
This contradicts our assumption \eqref{initialdata}.
   \item[iii)] If $\displaystyle n<\frac{1-\beta}{1+\beta}$, then
    $$p=\frac{n+1}{(n-1)_+}>-\frac{1}{\beta},$$
 i.e. $p\beta< -1$. In this case, we change the test function $\psi$ by $\displaystyle\psi(t)=\Psi\left(\frac{t}{K^{-1}T}\right)$ where $K\ge 1$ is independent of $T$.
Then 
$$\int_0^{K^{-1}T}(1+t)^{-\frac{(\beta+1)\,p}{p-1}}\,dt\leq C\,K^{-1+(\beta+1)p^\prime}\,T^{1-(\beta+1)p^\prime},\qquad\hbox{for all}\,\, T>0.$$ 
 From the subcritical case \eqref{10d11}, we can see easily that we have
 \begin{equation}\label{regularity3}
 u\in L^p((0,\infty);L^p(\mathbb{R}^n)).
 \end{equation}
  On the other hand, by applying H\"{o}lder's inequality in (\ref{10}) instead of Young's inequality, we get
\begin{eqnarray*}
\int_{\mathbb{R}^n}u_1(x)\phi_{T}(x)\,dx &\leq& C\,K^{2p^\prime-1}\left(\int_{\frac{K^{-1}T}{2}}^{K^{-1}T}\int_{\mathbb{R}^n}|u|^{p}\varphi\,dx\,dt\right)^{1/p}+\,C\,K^{-1}\left(\int_0^{K^{-1}T}\int_{|x|\geq T}|u|^{p}\varphi\,dx\,dt\right)^{1/p}\\
&{}&+ C\;K^{-1+(\beta+1)p^\prime}\,T^{-(\beta+1) p^\prime-\sigma p^\prime+1+n}+\,C\left(K^\sigma T^{-\sigma}+KT^{-1}\right)\,\int_{\mathbb{R}^n}|u_0(x)|\,dx.
\end{eqnarray*}
As $\sigma\geq1-\beta$ and $2>\sigma>\sigma+\beta$, we conclude that
 $$p=\frac{n+1}{n-1}\leq\frac{n+1}{n-\beta-\sigma}\Longrightarrow -(\beta+1) p^\prime-\sigma p^\prime+1+n\leq0,\quad\hbox{when}\,\,n\geq2,$$
 and
 $$-(\beta+1) p^\prime-\sigma p^\prime+1+n<0,\quad\hbox{when}\,\,n=1,$$
 so
\begin{eqnarray*}
\int_{\mathbb{R}^n}u_1(x)\phi_{T}(x)\,dx &\leq& C\,K^{2p^\prime-1}\left(\int_{\frac{K^{-1}T}{2}}^{K^{-1}T}\int_{\mathbb{R}^n}|u|^{p}\varphi\,dx\,dt\right)^{1/p}+\,C\,K^{-1}\left(\int_0^{K^{-1}T}\int_{|x|\geq T}|u|^{p}\varphi\,dx\,dt\right)^{1/p}\\&{}&+ C\;K^{-1+(\beta+1)p^\prime}+\,C\left(K^\sigma T^{-\sigma}+KT^{-1}\right)\,\int_{\mathbb{R}^n}|u_0(x)|\,dx.
\end{eqnarray*}
 Letting $T\longrightarrow\infty$ and taking into consideration \eqref{regularity3}, we get
 \begin{eqnarray*}
\int_{\mathbb{R}^n}u_1(x)\,dx &\leq& C\;K^{-1+(\beta+1)p^\prime}.
\end{eqnarray*}
 Letting $K\longrightarrow\infty$ and using $p\beta< -1\Longrightarrow -1+(\beta+1)p^\prime<0$, we infer that
$$\int_{\mathbb{R}^n}u_1(x)\,dx\leq 0.$$
This contradicts our assumption \eqref{initialdata}. 
 \end{enumerate}
\end{enumerate}


 \noindent {\bf II. Case: $\sigma<1-\beta$, for all $\beta<1$.}\\
In this case we take $\displaystyle d=\frac{1-\beta}{\sigma}>1$.\\
  
\noindent \underline{Subcritical case $\displaystyle p<\min\left\{\frac{n(1-\beta)+\sigma}{n(1-\beta)-\sigma};\frac{n}{n-\sigma}\right\}$.}\\
 Note that,
$$
\int_0^T(1+t)^{\displaystyle -\frac{(\beta+1)\,p}{p-1}}\,dt\leq C\,\left\{\begin{array}{ll}
T^{\displaystyle 1-\frac{(\beta+1)\,p}{p-1}}&\text{if}~\beta\,p<-1,\\
\ln T&\text{if}~\beta\,p=-1,\\
1&\text{if}~\beta\,p>-1.\\
\end{array}
\right.
$$
We have two cases to distinguish.
 \begin{enumerate}
 \item[a)] \underline{If $0\leq\beta<1$.} In this case, we have $\beta\,p>-1$ and so
 $$\int_0^T(1+t)^{-\frac{(\beta+1)\,p}{p-1}}\,dt\leq C.$$ 
 Therefore, (\ref{10}) implies
\begin{eqnarray*}
&{}&\int_{Q_T}|u|^p\varphi\,dx\,dt+\int_{\mathbb{R}^n}u_1(x)\phi_{T^d}(x)\,dx\\
&{}&\leq C\;T^{-2 p^\prime+1+\frac{n(1-\beta)}{\sigma}}+C\;T^{-\frac{2(1-\beta)}{\sigma}p^\prime+1+\frac{n(1-\beta)}{\sigma}}\\
&{}&\,\,\,\,+\,C\;T^{-(1-\beta)p^\prime+\frac{n(1-\beta)}{\sigma}}+
\,C\left(T^{-(1-\beta)}+T^{-1}\right)\,\int_{\mathbb{R}^n}|u_0(x)|\,dx.
\end{eqnarray*}
As 
 $$\sigma<1-\beta\Longrightarrow {-\frac{2(1-\beta)}{\sigma}p^\prime+1+\frac{n(1-\beta)}{\sigma}}< {-2 p^\prime+1+\frac{n(1-\beta)}{\sigma}},$$
 we conclude that
 \begin{eqnarray*}
&{}&\int_{Q_T}|u|^p\varphi\,dx\,dt+\int_{\mathbb{R}^n}u_1(x)\phi_{T^d}(x)\,dx\\
&{}&\leq C\;T^{-2 p^\prime+1+\frac{n(1-\beta)}{\sigma}}+\,C\;T^{(1-\beta)\left[-p^\prime+\frac{n}{\sigma}\right]}+
\,C\left(T^{-(1-\beta)}+T^{-1}\right)\,\int_{\mathbb{R}^n}|u_0(x)|\,dx,
\end{eqnarray*}
Note that,  we can easily see that
$${(1-\beta)\left[-p^\prime+\frac{n}{\sigma}\right]}<0\Longleftrightarrow p<\frac{n}{(n-\sigma)_+} \quad\hbox{and}\quad{-2 p^\prime+1+\frac{n(1-\beta)}{\sigma}}<0\Longleftrightarrow  p<\frac{n(1-\beta)+\sigma}{n(1-\beta)-\sigma}.$$
Note that $n(1-\beta)-\sigma\geq 1-\beta-\sigma>0$, for all $n\geq 1$. Letting $T\rightarrow\infty$, and using the Lebesgue
dominated convergence theorem together with $u_1\in L^1(\mathbb{R}^n)$, we conclude that
$$\int_{\mathbb{R}^n}u_1(x)\,dx\leq 0.$$
This contradicts our assumption \eqref{initialdata}.
 \item[b)] \underline{If $\beta<0$.} We have two cases for $n$.
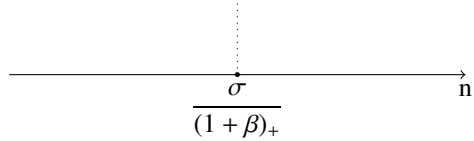
  \begin{center}
\begin{tikzpicture}[scale=1]
 \draw[->] (-3,0) -- (3,0) node[below]{n};
    \draw (0,0) node[below]{$\displaystyle\frac{\sigma}{(1+\beta)_+}$};
  \draw (0,0) node{{\tiny$\bullet$}};
   \draw[dotted](0,0)--(0,1);

        \end{tikzpicture}
 \end{center}
 Note that $n(1-\beta)-\sigma\geq1-\beta-\sigma>0$ and $\sigma<\beta+1<1\leq n\Longrightarrow n-\sigma>0$, for all $n\geq1$.
  \begin{enumerate}
 \item[i)] If $\displaystyle n\geq\frac{\sigma}{(1+\beta)_+}$, we have
  $$p<\min\left\{\frac{n(1-\beta)+\sigma}{n(1-\beta)-\sigma};\frac{n}{n-\sigma}\right\}=\frac{n}{n-\sigma}\leq-\frac{1}{\beta},$$
 i.e. $\beta\,p>-1$, which implies a contradiction by following the same calculations as in part a).
   \item[ii)] If $\displaystyle n<\frac{\sigma}{(1+\beta)_+}$, we get
   $$p<\min\left\{\frac{n(1-\beta)+\sigma}{n(1-\beta)-\sigma};\frac{n}{n-\sigma}\right\}=\frac{n(1-\beta)+\sigma}{n(1-\beta)-\sigma},$$
   and
   $$-\frac{1}{\beta}<\frac{n(1-\beta)+\sigma}{n(1-\beta)-\sigma},$$
   i.e.
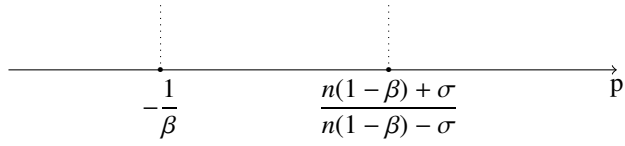
    \begin{center}
\begin{tikzpicture}[scale=1]
 \draw[->] (-4,0) -- (4,0) node[below]{p};
 \draw (1,0) node[below]{$\displaystyle\frac{n(1-\beta)+\sigma}{n(1-\beta)-\sigma}$};
  \draw (1,0) node{{\tiny$\bullet$}};
   \draw[dotted](1,0)--(1,1);
 \draw (-2,0) node[below]{$\displaystyle -\frac{1}{\beta}$};
  \draw (-2,0) node{{\tiny$\bullet$}};
   \draw[dotted](-2,0)--(-2,1);
        \end{tikzpicture}
 \end{center}

$\bullet$ If $p\leq -1/\beta$, i.e. $p\beta\geq -1$, then 
$$\int_0^T(1+t)^{-\frac{(\beta+1)\,p}{p-1}}\,dt\leq C\,\ln T,\qquad\hbox{for all}\,\, T\gg1.$$ 
Then, (\ref{10}) implies
\begin{eqnarray*}
&{}&\int_{Q_T}|u|^p\varphi\,dx\,dt+\int_{\mathbb{R}^n}u_1(x)\phi_{T^d}(x)\,dx\\
&{}&\leq C\;T^{-2 p^\prime+1+\frac{n(1-\beta)}{\sigma}}+C\;T^{-\frac{2(1-\beta)}{\sigma}p^\prime+1+\frac{n(1-\beta)}{\sigma}}\\
&{}&\,\,\,\,+\,C\;T^{-(1-\beta)p^\prime+\frac{n(1-\beta)}{\sigma}}\,\ln T+
\,C\left(T^{-(1-\beta)}+T^{-1}\right)\,\int_{\mathbb{R}^n}|u_0(x)|\,dx,
\end{eqnarray*}
 for all $T\gg1$. As 
 $$\sigma<1-\beta\Longrightarrow {-\frac{2(1-\beta)}{\sigma}p^\prime+1+\frac{n(1-\beta)}{\sigma}}< {-2 p^\prime+1+\frac{n(1-\beta)}{\sigma}},$$
 we conclude that
 \begin{eqnarray*}
&{}&\int_{Q_T}|u|^p\varphi\,dx\,dt+\int_{\mathbb{R}^n}u_1(x)\phi_{T^d}(x)\,dx\\
&{}&\leq C\;T^{-2 p^\prime+1+\frac{n(1-\beta)}{\sigma}}+\,C\;T^{(1-\beta)\left[-p^\prime+\frac{n}{\sigma}\right]}\,\ln T+
\,C\left(T^{-(1-\beta)}+T^{-1}\right)\,\int_{\mathbb{R}^n}|u_0(x)|\,dx,
\end{eqnarray*}
Note that,  we can easily see that
$${(1-\beta)\left[-p^\prime+\frac{n}{\sigma}\right]}<0\Longleftrightarrow p<\frac{n}{n-\sigma} \quad\hbox{and}\quad{-2 p^\prime+1+\frac{n(1-\beta)}{\sigma}}<0\Longleftrightarrow  p<\frac{n(1-\beta)+\sigma}{n(1-\beta)-\sigma}.$$
 Letting $T\rightarrow\infty$, using the fact that $\ln T\leq T^{\frac{{(1-\beta)\left(p^\prime-\frac{n}{\sigma}\right)}}{2}}$ (because $p<\frac{n}{n-\sigma}$) and the Lebesgue
dominated convergence theorem, we conclude that
$$\int_{\mathbb{R}^n}u_1(x)\,dx\leq 0.$$
This contradicts our assumption \eqref{initialdata}.\\

$\bullet$ If $p> -1/\beta$, i.e. $p\beta< -1$, then 
$$\int_0^T(1+t)^{-\frac{(\beta+1)\,p}{p-1}}\,dt\leq C\,T^{1-\frac{(\beta+1)\,p}{p-1}},\qquad\hbox{for all}\,\, T>0.$$ 
Then, (\ref{10}) implies
\begin{eqnarray*}
&{}&\int_{Q_T}|u|^p\varphi\,dx\,dt+\int_{\mathbb{R}^n}u_1(x)\phi_{T^d}(x)\,dx\\
&{}&\leq C\;T^{-2 p^\prime+1+\frac{n(1-\beta)}{\sigma}}+C\;T^{-\frac{2(1-\beta)}{\sigma}p^\prime+1+\frac{n(1-\beta)}{\sigma}}+
\,C\left(T^{-(1-\beta)}+T^{-1}\right)\,\int_{\mathbb{R}^n}|u_0(x)|\,dx.
\end{eqnarray*}
As 
 $$\sigma<1-\beta\Longrightarrow {-\frac{2(1-\beta)}{\sigma}p^\prime+1+\frac{n(1-\beta)}{\sigma}}< {-2 p^\prime+1+\frac{n(1-\beta)}{\sigma}},$$
 we conclude that
 \begin{eqnarray*}
&{}&\int_{Q_T}|u|^p\varphi\,dx\,dt+\int_{\mathbb{R}^n}u_1(x)\phi_{T^d}(x)\,dx\\
&{}&\leq C\;T^{-2 p^\prime+1+\frac{n(1-\beta)}{\sigma}}+\,C\left(T^{-(1-\beta)}+T^{-1}\right)\,\int_{\mathbb{R}^n}|u_0(x)|\,dx,
\end{eqnarray*}
Note that,  we can easily see that
$${-2 p^\prime+1+\frac{n(1-\beta)}{\sigma}}<0\Longleftrightarrow  p<\frac{n(1-\beta)+\sigma}{n(1-\beta)-\sigma}.$$
 Letting $T\rightarrow\infty$, using the Lebesgue
dominated convergence theorem, we conclude that
$$\int_{\mathbb{R}^n}u_1(x)\,dx\leq 0.$$
This contradicts our assumption \eqref{initialdata}.\\
 
 \end{enumerate}
\end{enumerate}

 \noindent \underline{Critical case: $\displaystyle p=\frac{n(1-\beta)+\sigma}{n(1-\beta)-\sigma}$ and $\sigma>n(1+\beta)$.}\\
Note that when $n(1+\beta)<\sigma$, we have
 $$p=\frac{n(1-\beta)+\sigma}{n(1-\beta)-\sigma}<\frac{n}{n-\sigma}\quad\hbox{i.e.}\quad  -\sigma p^\prime+n<0,$$
We have two cases to distinguish.
 \begin{enumerate}
 \item[a)] \underline{If $0\leq\beta<1$.} In this case, we have $\beta\,p>-1$ and so
 $$\int_0^T(1+t)^{-\frac{(\beta+1)\,p}{p-1}}\,dt\leq C.$$ 
 From the subcritical case \eqref{10d11}, we can see easily that we have
 \begin{equation}\label{regularity4}
 u\in L^p((0,\infty);L^p(\mathbb{R}^n)).
 \end{equation}
 On the other hand, by applying H\"{o}lder's inequality instead of Young's inequality, (\ref{10}) implies
\begin{eqnarray*}
\int_{\mathbb{R}^n}u_1(x)\phi_{T^d}(x)\,dx &\leq& C\left(\int_{\frac{T}{2}}^T\int_{\mathbb{R}^n}|u|^{p}\varphi\,dx\,dt\right)^{1/p}+\,C\left(\int_0^T\int_{|x|\geq T^d}|u|^{p}\varphi\,dx\,dt\right)^{1/p}\\
&{}&+ C\;T^{(1-\beta)\left[-p^\prime+\frac{n}{\sigma}\right]}+\,C\left(T^{-(1-\beta)}+T^{-1}\right)\,\int_{\mathbb{R}^n}|u_0(x)|\,dx.
\end{eqnarray*}
 Letting $T\longrightarrow\infty$ and taking into consideration \eqref{regularity4} we get
$$\int_{\mathbb{R}^n}u_1(x)\,dx\leq 0.$$
This contradicts our assumption \eqref{initialdata}.
 \item[b)] \underline{If $\beta<0$.} As $\displaystyle n(1+\beta)<\sigma\Longrightarrow n<\frac{\sigma}{(1+\beta)_+}$, so
   $$\displaystyle p=\frac{n(1-\beta)+\sigma}{n(1-\beta)-\sigma}>-\frac{1}{\beta},$$
     i.e. $p\beta< -1$. In this case, we change the test function $\psi$ by $\displaystyle \psi(t)=\Psi\left(\frac{t}{K^{-1}T}\right)$ where $K\ge 1$ is independent of $T$.
Then 
     $$\displaystyle\int_0^{K^{-1}T}(1+t)^{-\frac{(\beta+1)\,p}{p-1}}\,dt\leq C\,K^{-1+(\beta+1)p^\prime}\,T^{1-(\beta+1)p^\prime},\qquad\hbox{for all}\,\, T>0.$$ 
     From the subcritical case \eqref{10d11}, we can see easily that we have
 \begin{equation}\label{regularity5}
 u\in L^p((0,\infty);L^p(\mathbb{R}^n)).
 \end{equation}
  On the other hand, by applying H\"{o}lder's inequality in  (\ref{10}) instead of Young's inequality and using $\displaystyle p=\frac{n(1-\beta)+\sigma}{n(1-\beta)-\sigma}$, (\ref{10}) implies
\begin{eqnarray*}
\displaystyle
&{}&\displaystyle\int_{Q_T}|u|^p\varphi\,dx\,dt+\int_{\mathbb{R}^n}u_1(x)\phi_{T^d}(x)\,dx\\
&{}&\displaystyle\leq C\,K^{2p^\prime-1}\left(\int_{\frac{K^{-1}T}{2}}^{K^{-1}T}\int_{\mathbb{R}^n}|u|^{p}\varphi\,dx\,dt\right)^{1/p}+\,C\,K^{-1}\left(\int_0^{K^{-1}T}\int_{|x|\geq T^d}|u|^{p}\varphi\,dx\,dt\right)^{1/p}\\
&{}&\displaystyle\,\,+C\;K^{-1+(\beta+1)p^\prime}+
\,C\left(K^{1-\beta}T^{-(1-\beta)}+KT^{-1}\right)\,\int_{\mathbb{R}^n}|u_0(x)|\,dx.
\end{eqnarray*}
 Letting $T\longrightarrow\infty$ and taking into consideration \eqref{regularity5}, we get
 \begin{eqnarray*}
\int_{\mathbb{R}^n}u_1(x)\,dx &\leq& C\;K^{-1+(\beta+1)p^\prime}.
\end{eqnarray*}
 Letting $K\longrightarrow\infty$ and using $p\beta< -1\Longrightarrow -1+(\beta+1)p^\prime<0$, we infer that
$$\int_{\mathbb{R}^n}u_1(x)\,dx\leq 0.$$
This contradicts our assumption \eqref{initialdata}.\\

\hfill$\square$
\end{enumerate}

\section*{Acknowledgments} 

M.K. gratefully acknowledges the financial support by KUST through the grant FSU-2021-015.
A.L. gratefully acknowledges the financial support by KUST through the grant FSU-2021-027.



\end{document}